\documentclass[12pt]{article}
\usepackage[T1]{fontenc}
\usepackage{color}

\usepackage[english,francais]{babel}
\usepackage{authblk}

\usepackage{amssymb,url,xspace,amsmath}



\newcommand{\BibTeX}{{\scshape Bib}\kern-.08em\TeX}
\newcommand{\T}{\S\kern .15em\relax }
\newcommand{\AMS}{$\mathcal{A}$\kern-.1667em\lower.5ex\hbox
        {$\mathcal{M}$}\kern-.125em$\mathcal{S}$}

\title{D\'ecomposition effective de Jordan-Chevalley}
\author[1]{Danielle Couty\thanks{danielle.couty@iut-tarbes.fr}}
\author[2,3]{Jean Esterle\thanks{esterle@math.u-bordeaux1.fr,  j.esterle@estia.fr}}
\author[4]{Rachid Zarouf\thanks{rzarouf@cmi.univ-mrs.fr}}
\affil[1]{IMT Toulouse, UMR 5219, 118 route de Narbonne, 31062, Toulouse Cedex 9}
\affil[2]{Universit\'e de Bordeaux, IMB, UMR 5251, 351 Cours de la Lib\'eration, 33405 Talence Cedex}
\affil[3]{ESTIA, Technopole Izarbel, 64200, Bidart}
\affil[4]{CMI-LATP, UMR 6632, Universit\'e de Provence, 39, rue F.-Joliot-Curie, 13453 Marseille cedex 13}

\date {16, juin 2011}

\begin{document}

\maketitle

\section{Introduction}

On note ${\mathcal{M}}_{n}(k)$ l'alg\`ebre des matrices $n\times n$
\`a coefficients dans un corps $k,$ et on  note ${\mathcal{G}}L_{n}(k)$
le groupe des \'el\'ements inversibles de ${\mathcal{M}}_{n}(k).$ Une
matrice $V\in{\mathcal{M}}_{n}(k)$ est dite unipotente si $I_{n}-V$
est nilpotente, $I_{n}$ d\'esignant la matrice unit\'e. Le th\'eor\`eme de
d\'ecomposition de Jordan, sous sa forme additive, montre que toute
matrice $U\in{\mathcal{M}}_{n}(k)$ s'\'ecrit de mani\`ere unique sous
la forme $U=D+N,$ où $D\in{\mathcal{M}}_{n}(\tilde{k})$ est diagonalisable
sur le corps de d\'ecomposition $\tilde{k}$ du polynôme caract\'eristique
$p_{U}$ de $U,$ et où $N$ est nilpotente et commute avec $D.$
De plus il existe $p\in k[x]$ tel que $D=p(U),$ $N=U-p(U).$

Si $k$ est un corps parfait, c'est-\`a-dire si tout polynôme irr\'eductible
$p\in k[x]$ est \`a racines simples dans son corps de d\'ecomposition,
alors $D$ est $N$ appartiennent \`a ${\mathcal{M}}_{n}(k).$ Cette
propri\'et\'e reste v\'erifi\'ee sans hypoth\`ese particuli\`ere sur $k$ quand
tout facteur irr\'eductible de $p_{U}$ est \`a racines simples dans son
corps de d\'ecomposition.

En posant $V=I+D^{-1}N,$ on d\'eduit imm\'ediatement de la forme additive
de la d\'ecomposition de Jordan que si $k$ est parfait alors toute
matrice $U\in{\mathcal{G}}L_{n}(k)$ s'\'ecrit de mani\`ere unique sous
la forme

\begin{equation}
U=DV,\end{equation}
 avec $D\in{\mathcal{G}}L_{n}(k)$ diagonalisable sur le corps de
d\'ecomposition de $p_{U},$ et $V$ unipotente commutant avec $D.$

Le but de cet article est d'attirer l'attention sur le fait que ces
d\'ecompositions de Jordan sont \textit{effectivement calculables} \`a
partir du polynôme caract\'eristique $p_{U}$ de $U.$ La m\'ethode, qui
date du d\'ebut des ann\'ees 50, est due \`a Claude Chevalley. Elle est
directement inspir\'ee de la m\'ethode de Newton, ou \textquotedbl{}m\'ethode
de la tangente\textquotedbl{}, qui fournit un sch\'ema d'approximation
d'une racine de l'\'equation $f(x)=0$. La m\'ethode de Newton consiste
\`a choisir convenablement $x_{0}$ et \`a poser, pour $m\geq0$, \[
x_{m+1}=x_{m}-\frac{f(x_{m})}{f'(x_{m})}.\]
 Des hypoth\`eses classiques permettent de garantir que la suite $(x_{m})_{m\geq0}$
converge vers une solution $x$ de l'\'equation $f(x)=0$. L'id\'ee qui
sous-tend l'algorithme de Chevalley est de trouver de cette façon
une solution non triviale $D$ de l'\'equation $\tilde{p}(D)=0$, où
$\tilde{p}$ est le produit des facteurs irr\'eductibles distincts,
\`a racines distinctes dans leur corps de d\'ecomposition, d'un polynôme
annulateur $p$ de la matrice $U.$ Il n'y a aucun probl\`eme de convergence
ici, car la suite obtenue est stationnaire: si $p=p_{1}^{n_{1}}\dots p_{m}^{n_{m}}$
est la d\'ecomposition en produit de facteurs irr\'eductibles distincts
de $p,$ et si on pose

\[
D_{0}=U,\ D_{n+1}=D_{n}-\tilde{p}(D_{n})[\tilde{p}'(D_{n})]^{-1},\]
 alors la suite $(D_{n})_{n\ge0}$ est bien d\'efinie, et $D_{n}=D$
pour $2^{n}\ge$max$_{1\le j\le m}n_{j}$, (voir Section 2,  Th\'eor\`eme 1). Le fait que $\tilde{p}'(D_n)$
est inversible r\'esulte du fait que $\tilde{p}'$ est premier avec $\tilde{p},$ donc
$p'(D_{n})^{-1}=q(D_{n}),$ où $q\in k[x]$ est un inverse de $\tilde{p}'$
mod $p$, car on v\'erifie que $p(D_{n})=0$ pour $n\ge0$ (en fait
on peut prendre un inverse de $\tilde{p}'$ mod $p_{1}^{n_{1}-1}\dots p_{m}^{n_{m}-1}).$
On trouvera la version originale de cet algorithme (th\'eor\`eme 7 page
71) dans l'ouvrage publi\'e en 1951 par C.Chevalley \cite{Che51}.

Cette m\'ethode est bien connue dans les pr\'eparations \`a l'agr\'egation
\cite{Fe}, \cite{Pic}, \cite{FrMa}, et est par exemple propos\'ee
en exercice p.62-63 de l'ouvrage r\'ecent de A.Boyer et J.Risler \cite{BoRi}, mais elle ne semble pas avoir eu jusqu'ici la diffusion qu'elle
m\'erite dans les enseignements au niveau L2-L3. Beaucoup de coll\`egues
sont \'etonn\'es d'apprendre que le calcul de la d\'ecomposition de Jordan-Chevalley,
souvent appel\'ee d\'ecomposition de Dunford, ne n\'ecessite pas la connaissance
des valeurs propres de la matrice consid\'er\'ee. 

La forme multiplicative de la d\'ecomposition de Jordan-Chevalley joue
un rôle important dans la th\'eorie des groupes alg\'ebriques, c'est \`a
dire des sous-groupes $G$ de ${\mathcal{G}}L_{n}(k)$ de la forme
\[
G=\left\{ U=(u_{i,j})_{\stackrel{1\le i\le n}{_{1\le j\le n}}}\ |\ \phi_{\lambda}[u_{1,1},\, u_{1,2},\dots,\, u_{n,n-1},\, u_{n,n}]=0,\lambda\in\Lambda_{G}\right\} ,\]
 où ($\Lambda_G$ est une famille d'indices et) $(\phi_{\lambda})_{\lambda\in\Lambda_{G}}$ est une famille de
polynômes en $n^{2}$ variables \`a coefficients dans $k.$ En effet
si $k$ est parfait, et si $G\subset{\mathcal{G}}L_{n}(k)$ est un
groupe alg\'ebrique, alors la partie diagonalisable $D$ et la partie
unipotente $V$ de la d\'ecomposition de Jordan-Chevalley $U=DV$ d'un
\'el\'ement $U$ de $G$, appartiennent \`a $G.$ Ce fait a \'et\'e largement
utilis\'e par le premier auteur dans sa th\`ese \cite{Cou}. Outre l'ouvrage
fondateur de Chevalley \cite{Che51} d\'ej\`a mentionn\'e plus haut, nous
renvoyons aux monographies de Borel \cite{Bor56} et Humphreys \cite{Hum}
pour une pr\'esentation g\'en\'erale de la th\'eorie des groupes alg\'ebriques.

Nous r\'esumons maintenant bri\`evement la suite de l'article. Dans la
partie II nous donnons une version tr\`es g\'en\'erale de la d\'ecomposition
de Jordan-Chevalley en nous pla\c cant de m\^eme que dans \cite{Duc} dans le cadre d'une alg\`ebre $A$ sur un corps $k.$ 
En reprenant une terminologie introduite dans \cite{Bou}, Chapitre VII pour les endomorphismes, on dit que $u\in A$ est absolument semi-simple (resp. semi-simple, resp. s\'eparable) s'il poss\`ede
un polynôme annulateur premier \`a son polynôme d\'eriv\'e (resp. dont les facteurs irr\'eductibles sont distincts, resp. dont les
facteurs irr\'eductibles sont premiers \`a leur polynôme d\'eriv\'e). On montre alors que tout \'el\'ement s\'eparable
$u$ de $A$ s'\'ecrit de mani\`ere unique sous la forme $u=d+n,$ avec $d\in A$
absolument semi-simple et $n\in A$ nilpotent. La d\'emonstration est essentiellement
la d\'emonstration originale de Chevalley, même si pour la preuve de
l'unicit\'e nous donnons une d\'emonstration directe qui \'evite le recours
habituel au corps de d\'ecomposition du polynôme minimal de $u.$

 On conclut cette section en rappelant la caract\'erisation des corps parfaits et en donnant des exemples d'endomorphismes semi-simples  qui ne sont pas absolument semi-simples. De tels endomorphismes n'admettent pas de d\'ecomposition de Jordan, comme l'a observ\'e Bourbaki dans \cite{Bou}, Chapitre VII.

A la section III nous illustrons le caract\`ere effectif de l'algorithme
de Chevalley sur un exemple : quand $k$ est de caract\'eristique nulle
le polynôme $\tilde{p}$ utilis\'e dans l'algorithme est \'egal \`a $\frac{p}{pgcd(p,p')},$
et dans le cas d'une matrice $U$ on peut prendre $p=p_{U}$ où $p_{U}$
d\'esigne le polynôme caract\'eristique de $U.$ On est donc ramen\'e via
l'algorithme d'Euclide \'etendu \`a une suite de divisions euclidiennes
pr\'ec\'ed\'ee d'un calcul de d\'eterminant. Avec l'aide de Maple nous montrons
ainsi comment calculer la d\'ecomposition de Jordan-Chevalley d'une
matrice $U\in{\mathcal M}_{15}(\mathbb{R})$ dont les 5 racines, de multiplicit\'e
3, ne sont pas calculables.

A la section IV nous discutons une autre m\'ethode, utilis\'ee par A.
Borel dans \cite{Bor56} et bas\'ee sur un syst\`eme ${\mathcal{S}}$
d'\'equations de congruence, qui permet de montrer l'existence de la
d\'ecomposition de Jordan (en utilisant un peu de th\'eorie de Galois
pour montrer que le polynôme permettant de calculer le terme absolument semi-simple
$d$ de la d\'ecomposition est \`a coefficients dans $k$) et de la calculer
explicitement quand on connaît les racines d'un polynôme annulateur.
La partie absolument semi-simple $d$ d'un \'el\'ement s\'eparable $u$ d'une $k$-alg\`ebre
$A$ s'\'ecrit alors sous la forme $d=q(u),$ où $q$ est une solution
quelconque de ${\mathcal{S}}.$ De même que dans \cite{Pic} nous
montrons que l'algorithme de Newton permet de calculer une solution
$q$ du syst\`eme ${\mathcal{S}}$ sans faire intervenir les racines.
Il n'y a pas de nouveau calcul \`a faire puisque $q$ s'obtient en appliquant
directement l'algorithme de la section II \`a $\pi(x),$ où $p$ d\'esigne
un polynôme annulateur de $u$ dont les facteurs irr\'eductibles sont
\`a d\'eriv\'ee non nulle et où $\pi:k[x]\to k[x]/pk[x]$ d\'esigne la surjection
canonique. En utilisant Maple nous donnons la valeur du polynôme $p$
de degr\'e 14 tel que la partie diagonalisable $D$ de la matrice $U$
\'etudi\'ee \`a la section III soit \'egale \`a $p(U).$

Nous indiquons d\`es maintenant pourquoi nous avons pr\'ef\'er\'e appeler
la d\'ecomposition discut\'ee dans cet article \textquotedbl{}d\'ecomposition
de Jordan-Chevalley\textquotedbl{} plutôt que \textquotedbl{}d\'ecomposition
de Dunford\textquotedbl{}. On trouve bien dans un article \'ecrit en
1954 par N. Dunford   \cite{Dun} une d\'ecomposition de type Jordan en
th\'eorie spectrale \footnote{Les r\'esultats de 1954 sont largement diffus\'es \`a partir de 1971, date de la parution de l'ouvrage de Dunford et Schwartz\emph{ Linear operators} \cite{DunSch}. Les auteurs pr\'esentent les r\'esultats de 1954 de Dunford au tome III, paragraphe "The Canonical Reduction of a Spectral Operator", pages 1937-1941, (voir aussi la note page 2096) et les g\'en\'eralisent (\cite{DunSch}, Th\'eor\`eme 8 p.2252).}, mais cet article est post\'erieur \`a l'ouvrage de
Chevalley cit\'e plus haut. Ceci justifie pour nous la terminologie
de \textquotedbl{}d\'ecomposition de Jordan-Chevalley\textquotedbl{},
ce que nous d\'evelopperons davantage dans la partie V en \'etudiant quelques
\'etapes historiques du \textquotedbl{}th\'eor\`eme de Jordan\textquotedbl{}.
Cette terminologie est d'ailleurs la terminologie utilis\'ee par James
E. Humphreys dans les deux ouvrages où il pr\'esente cette d\'ecomposition,
voir \cite{Hum72}, page 17 et \cite{Hum}, page 18. Nous compl\'etons cette pr\'esentation historique \`a la section VI en mentionnant deux aspects de la d\'ecomposition multiplicative de Jordan dans les groupes de Lie semi-simples: la d\'ecomposition d'un \'el\'ement d'un groupe de Lie r\'eel semi-simple et connexe en produit d'un \'el\'ement "elliptique" $e$, d'un \'el\'ement "hyperbolique" $h$ et d'un \'el\'ement unipotent $u$ qui commutent entre eux, et un r\'esultat tr\`es r\'ecent de
Venkataramana \cite{Ve} qui montre que si $G$ est un groupe de Lie lin\'eaire semi-simple r\'eel et si $g=su$ est la d\'ecomposition de Jordan d'un \'el\'ement $g$ d'un r\'eseau $\Gamma$ de $G$ alors il existe $m\ge 1$ tels que
$s^m\in \Gamma$ et $u^m \in \Gamma.$

Enfin dans la section VII nous discutons l'int\'erêt d'utiliser la d\'ecomposition
de Jordan-Chevalley pour des enseignements \`a un niveau relativement
\'el\'ementaire, en nous appuyant sur l'exp\'erience acquise \`a l'Ecole d'Ing\'enieurs
ESTIA de Bidart, dont les promotions (150 \'el\`eves de 1e ann\'ee en 2010-2011)
sont compos\'ees pour 2/3 d'\'el\`eves issus de diverses classes pr\'eparatoires,
l'autre tiers \'etant majoritairement form\'e de titulaires d'un DUT.

\section{L'algorithme de Newton pour la d\'ecomposition de Jordan-Chevalley}

Soit $k$ un corps, soit $k[x]$ l'anneau des polynômes \`a coefficients
dans $k.$ On dit qu'un polynôme non constant $p\in k[x]$ est \textit{s\'eparable}
s'il est premier \`a son polynôme d\'eriv\'e $p'$ (c'est \'evidemment le
cas si $p$ est irr\'eductible et si $p'\neq0).$ Il r\'esulte imm\'ediatement
de l'identit\'e de Bezout que si $q$ est un diviseur non constant d'un
polynôme s\'eparable $p$ alors $q$ est s\'eparable. Si $p$ et $q$
sont s\'eparables et premiers entre eux, alors $p$ est premier \`a $p'q,$
donc \`a $p'q+q'p=(pq)'.$ De même $q$ est premier \`a $(pq)',$ et par
cons\'equent $pq$ est premier avec $(pq)'$ donc $pq$ est s\'eparable.
On en d\'eduit imm\'ediatement que le produit d'une famille finie de polynômes
s\'eparables premiers entre eux deux \`a deux est s\'eparable.

Soit $A$ une $k$-alg\`ebre g\'en\'erale unitaire d'unit\'e not\'ee $1_{A}.$ On pose
$[u,\, v]=uv-vu$ pour $u,v\in A,$ $p(u)=p_{0}1_{A}+p_{1}u+...+p_{n}u^{n}$
pour $p=p_{0}+p_{1}x+\dots+p_{n}x^{n}\in k[x]$ et $u\in A,$ et on
dit que $p$ est un polynôme annulateur de $u$ si $p(u)=0.$ Par convention, les facteurs irr\'eductibles d'un polyn\^ome $p\in k[x]$ sont suppos\'es unitaires.
Les notions suivantes sont une extension aux \'el\'ements d'une $k$-alg\`ebre g\'en\'erales de notions introduites dans \cite{Bou}, Chapitre VII pour les endomorphismes d'un $k$-espace vectoriel.

\textbf{Definition 1.}  \textit{On dit que $u\in A$ est }

\textit{( i) absolument semi-simple
si $u$ poss\`ede un polynôme annulateur s\'eparable,}

\textit{(ii) semi-simple si $u$ poss\`ede un polyn\^ome annulateur dont tous les facteurs irr\'eductibles sont distincts,}

\textit{iii) s\'eparable si $u$ poss\`ede un
polynôme annulateur $p\in k[x]$ dont tous les facteurs irr\'eductibles
sont s\'eparables.}

 \textit {On dit qu'un couple $(d,\, s)$ d'\'el\'ements de $A$, avec
$d$ absolument semi-simple, $s$ nilpotent est une d\'ecomposition de Jordan-Chevalley
de $u\in A$ quand $u=d+s$ et $[d,\, s]=0.$}

Soit $u$ un \'el\'ement de $A$ poss\'edant un polynôme annulateur. Comme
tout diviseur non constant d'un polynôme s\'eparable est s\'eparable,
on voit que $u$ est absolument semi-simple si et seulement si son polynôme minimal
est s\'eparable, et que $u$ est s\'eparable si et seulement si tous les
facteurs irr\'eductibles de son polynôme minimal sont s\'eparables.

Notons que si $u\in A$ admet une d\'ecomposition de Jordan-Chevalley
$(d,\, s)$ alors $u$ est s\'eparable. En effet dans ce cas soit $p$
le polynôme minimal de $d.$ On a $p(u)=p(d)+sv=sv,$ avec $v\in A,$
$[s,\, v]=0,$ et $sv$ est nilpotent. Donc il existe un entier $m\ge1$
tel que $p(u)^{m}=0.$ Comme tout diviseur irr\'eductible de $p^{m}$
divise $p,$ $u$ est s\'eparable.

\smallskip{}

Il est bien connu que si $u$ et $v$ sont absolument semi-simples, et si $[u,\, v]=0,$
alors $u+v$ est absolument semi-simple. Nous donnons ici un r\'esultat qui est
suffisant pour montrer l'unicit\'e de la d\'ecomposition de Jordan, et
dont la d\'emonstration ne fait pas appel \`a la notion de clôture alg\'ebrique.

\textbf{Proposition 1.} \textit{Soient $u$ et $v$ deux \'el\'ements
absolument semi-simples d'une alg\`ebre unitaire $A$ sur un corps $k,$ et soit
$k[u,\, v]$ la sous-alg\`ebre unitaire de $A$ engendr\'ee par $u$ et
$v.$ Si $[u,\, v]=0,$ alors $k[u,\, v]$ ne contient aucun \'el\'ement
nilpotent non nul.}

Preuve. Soit $p\in k[x]$ le polynôme minimal de $u.$ Comme $u$
est absolument semi-simple, il existe une famille $p_{1},\dots,\, p_{s}$ de polynômes
unitaires irr\'eductibles distincts sur $k$ tels que $p=p_{1}\dots p_{s}.$
Posons $P_{i}=\prod_{j\neq i}p_{j}.$ D'apr\`es le th\'eor\`eme de Bezout
il existe $\phi_{1},\dots,\,\phi_{s}\in k[x]$ tels que $\sum_{j=1}^{s}\phi_{i}P_{i}=1_{k}.$
Posons $e_{i}=\phi_{i}(u)P_{i}(u)$, $u_{i}=e_{i}u$ et $A_{i}=e_{i}k[u,\, v].$
Alors $\sum_{i=1}^{s}e_{i}=1_{A}$, $e_{i}\neq0,$ $e_{i}e_{j}=0$ pour $i\neq j,$ ce qui implique $e_{i}^{2}=e_{i}\neq0$ pour tout $i;$   $A_{i}$
est unitaire d'unit\'e $e_{i}$ et $k[u,\, v]=\oplus_{1\le i\le s}A_{i}.$
Comme $e_{i}p_{i}(e_{i}u)=e_{i}p_{i}(u)=0,$ le polynôme minimal de
$u_{i}$ consid\'er\'e comme \'el\'ement de la $k$-alg\`ebre unitaire $A_{i}$
divise $p_{i},$ donc est \'egal \`a $p_{i}.$ Donc $k[u_{i}]\subset A_{i}$
est un corps isomorphe \`a l'alg\`ebre quotient $k[x]/p_{i}k[x].$

Consid\'erons maintenant $A_{i}$ comme une $k[u_{i}]$-alg\`ebre unitaire.
Alors $v_{i}:=e_{i}v$ est un \'el\'ement absolument semi-simple de $A_{i}$ et en
raisonnant de même que ci-dessus on construit une famille $(e_{i,j})_{1\le j\le\alpha_{i}}$
d'\'el\'ements non nuls de $A_{i},$ avec $\alpha_{i}\ge1,$ v\'erifiant
$e_{i,j}e_{i,j'}=0$ pour $j'\neq j$ et $\sum_{j=1}^{\alpha_{i}}e_{i,j}=e_{i}$
et tels que $A_{i,j}:=e_{i,j}A_{i}=e_{i,j}k[u,v]=e_{i,j}k[u_{i}][v]=e_{i,j}k[u_{i}][e_{i,j}v]$
soit un corps d'unit\'e $e_{i,j}.$

Soit maintenant $w$ un \'el\'ement nilpotent de $k[u,\, v].$ Alors $e_{i,j}w$
est nilpotent, donc $e_{i,j}w=0$ pour $1\le j\le s,1\le j\le\alpha_{i}$
et $w=1_{k}w=\sum\limits _{i=1}^{n}\left(\sum\limits _{j=1}^{\alpha_{i}}e_{i,j}w\right)=0.$
$\square$

On va maintenant d\'emontrer une version tr\`es g\'en\'erale du th\'eor\`eme de
d\'ecomposition de Jordan-Chevalley, en utilisant un algorithme, directement
inspir\'e par la m\'ethode de Newton, introduit par C. Chevalley. 

\begin{flushleft}
\textbf{Th\'eor\`eme 1. }\textit{Soient $k$ un corps, $A$ une $k$-alg\`ebre,
$u$ un \'el\'ement s\'eparable de $A,$ et $p=p_{1}^{m_{1}}\dots p_{r}^{m_{r}}\in k[x]$
un polynôme annulateur unitaire de $u$ dont les facteurs irr\'eductibles
unitaires $p_{1},\dots,\, p_{r}$ sont s\'eparables. On pose $m=$$\max_{1\le j\le r}m_{j},$
$\tilde{p}=p_{1}\dots p_{r},$ $\overline{p}=p/\tilde{p},$ de sorte
que $\tilde{p}$ est s\'eparable. Il existe $q\in k[x]$ tel que $q\tilde{p}'\equiv1\ {\rm mod}\ \overline{p}.$
On d\'efinit par r\'ecurrence une suite $(d_{n})_{n\geq0}$ d'\'el\'ements
de $A$ en posant\[
\left\{ \begin{array}{c}
d_{0}=u\\
d_{n+1}=d_{n}-\tilde{p}(d_{n})q(d_{n})\; pour\; n\geq0.\end{array}\right.\]
 Soit $N$ le plus petit entier tel que $2^{N}\ge m.$ Alors $\tilde{p}(d_{N})=0,$
$d_{N}$ est absolument semi-simple, $u-d_{N}$ est nilpotent, $[d_{N},\, u-d_{N}]=0$
et $(d_{N},\, u-d_{N})$ est l'unique d\'ecomposition de Jordan-Chevalley
de $u.$}
\par\end{flushleft}

\begin{flushleft}
Preuve. Comme $\tilde{p}'$ est premier avec $\tilde{p},$ il est
premier avec $p_{j}$ pour $1\le j\le r,$ donc il est premier avec
$p$ et $\overline{p},$ ce qui permet de construire $q.$
\par\end{flushleft}

\begin{flushleft}
Soit $B:=\{h(u)\}_{h\in k[x]}.$ Il est clair que $B$ est commutative,
et que $d_{n}\in B$ pour $n\ge0.$
\par\end{flushleft}

\begin{flushleft}
Observons tout d'abord que si $f\in k[x],$ on a, pour $v,\, w\in B,$
\par\end{flushleft}

\begin{flushleft}
\begin{equation}
f(v+w)-f(v)\in wB,\ f(v+w)-f(v)-wf'(v)\in w^{2}B.\end{equation}

\par\end{flushleft}

\begin{flushleft}
Si $v\in B$ v\'erifie $p(v)=0,$ et si $w\in\tilde{p}(v)B,$ alors
$p_{j}(v+w)-p_{j}(v)\in wB\subset p_{j}(v)B,$ donc $p_{j}(v+w)\in p_{j}(v)B$
pour $1\le j\le r,$ ce qui montre que $p(v+w)=0.$
\par\end{flushleft}

\begin{flushleft}
Comme $d_{0}=u,$ et comme $p$ est un polynôme annulateur de $u,$
une r\'ecurrence imm\'ediate montre que $p(d_{n})=0$ pour $n\ge1.$
\par\end{flushleft}

\begin{flushleft}
On a alors,
\par\end{flushleft}

\begin{flushleft}
\[
\tilde{p}'(d_{n})q(d_{n})-1_{A}\in\overline{p}(d_{n})B,\: q(d_{n})-\tilde{p}'(d_{n})^{-1}\in\overline{p}(d_{n})B,\ \tilde{p}(d_{n})q(d_{n})=\tilde{p}(d_{n})\tilde{p}'(d_{n})^{-1}.\]
 On obtient, pour $n\ge0,$
\par\end{flushleft}

\begin{flushleft}
\begin{equation}
d_{n+1}=d_{n}-\tilde{p}(d_{n})\tilde{p}'(d_{n})^{-1},\end{equation}

\par\end{flushleft}

\begin{flushleft}
ce qui montre que la suite $(d_{n})_{n\ge0}$ est d\'efinie par l'algorithme
de Newton associ\'e \`a $\tilde{p}.$
\par\end{flushleft}

\begin{flushleft}
On a $d_{0}=u$ et, pour $n\ge1,$
\par\end{flushleft}

\begin{flushleft}
\[
u-d_{n}=\sum_{j=0}^{n-1}\tilde{p}(d_{j})q(d_{j}),\]

\par\end{flushleft}

\begin{flushleft}
ce qui prouve que $(u-d_{n})^{m}=0$ pour $n\ge0.$ De plus, il
r\'esulte de la formule (2) que $\tilde{p}(d_{n+1})\in\tilde{p}(d_{n})^{2}B.$
Par r\'ecurrence, on en deduit que pour tout $n\in\mathbb{N},$ \[
\tilde{p}(d_{n})\in\tilde{p}(u)^{2^{n}}B.\]
 En particulier, $\tilde{p}(d_{N})=0.$ Donc $d_{N}$ est absolument semi-simple,
et comme $u-d_{N}$ est nilpotent on voit que $(d_{N},\, u-d_{N})$
est une d\'ecomposition de Jordan-Chevalley de $u.$ Soit $(d,\, s)$
une autre d\'ecomposition de Jordan-Chevalley de $u.$ Comme $d_{N}\in B,$
$[d,\, d_{N}]=[s,\, u-d_{N}]=0$ et $d-d_{N}\in k[d,\, d_{N}]$ est
nilpotent. Il r\'esulte alors de la proposition 1 que $d=d_{N},$ et
la d\'ecomposition de Jordan de $u$ est unique. $\square$
\par\end{flushleft}

\bigskip{}

Rappelons qu'on dit qu'un corps $k$ est parfait si tout polynôme
irr\'eductible $p\in k[x]$ est s\'eparable. On a alors le corollaire
suivant. \smallskip{}

\textbf{Corollaire 1.} \textit{Soit $k$ un corps parfait, soit $A$
une $k$-alg\`ebre et soit $u\in A$. Si $u$ admet un polynôme annulateur,
alors $u$ admet une unique d\'ecomposition de Jordan-Chevalley $(d,\, s)$
sur $A$, et il existe $p\in k[x]$ tel que $d=p(u),\, s=u-p(u).$}

\smallskip{}

Il est clair que tout corps de caract\'eristique nulle est parfait,
et un corps $k$ de caract\'eristique $m\ge2$ est parfait si l'\'equation
$y^{m}=a$ admet une solution dans $k$ pour tout $a\in k.$ En effet
dans ce cas soit $p\in k[x]$ un polynôme non constant tel que $p'=0.$
Alors il existe $s\ge1$ et $\lambda_{0},\dots,\,\lambda_{s}\in k$
tels que $p=\sum\limits _{j=0}^{s}\lambda_{j}x^{mj}=\left(\sum\limits _{j=0}^{s}\mu_{j}x^{j}\right)^{m},$
où $\mu_{j}$ est une solution dans $k$ de l'\'equation $z^{m}=\lambda_{j}$
pour $0\le j\le s,$ donc $p$ n'est pas irr\'eductible, ce qui montre
que $k$ est parfait. R\'eciproquement soit $k$ un corps de caract\'eristique
$m\ge1$ tel qu'il existe $a\in k$ v\'erifiant $y^{m}\neq a$ pour
tout $y\in k.$ Posons $p=x^{m}-a,$ soit $\tilde{k}$ un corps de
rupture de $p$ et soit $b$ une racine de $p$ dans $\tilde{k}.$
Alors $p=x^{m}-a=x^{m}-b^{m}=(x-b)^{m},$ donc les seuls diviseurs
unitaires non constants de $p$ dans $k[x]$ distincts de $p$ sont
de la forme $(x-b)^{n}=x^{n}-nbx^{n-1}+\sum\limits _{j=2}^{n}C_{n}^{j}(-n)^{j}x^{n-j},$
avec $1\le n<m.$ Comme dans ce cas $n$ est inversible modulo $m,$
aucun de ces polynômes n'appartient \`a $k[x]$ et $p$ est un polynôme
irr\'eductible de d\'eriv\'ee nulle, ce qui montre que $k$ n'est pas parfait.

Le corollaire 1 n'est plus valable si $k$ n'est pas parfait. En effet soit $k$ un corps non parfait, de caract\'eristique $m\neq 0,$ soit $a \in k$ tel que l'\'equation $x^m=a$ n'admette aucune solution dans $k$, et soit $E$ un $k$-espace vectoriel de dimension $m.$ Soit $\{e_1,...,e_m\}$ une base de $E,$ et soit $\theta \in {\mathcal L}(E)$ l'endomorphisme d\'efini par les formules $\theta(e_{j})=e_{j+1}$ pour $1\le j\le m-1$ et $\theta(e_{m})=ae_{1}.$
Alors $\theta^{m}=a1_{{\mathcal{L}}(E)}$ et $x^{m}-a$ divise tout
polynôme annulateur de $\theta,$ ce qui prouve que $\theta$ est
un \'el\'ement non s\'eparable de la $k$-alg\`ebre ${\mathcal{L}}(E).$ Donc $\theta$ n'admet pas de d\'ecomposition de Jordan. En fait cet endomorphisme $\theta$ est un exemple d'endomorphisme semi-simple qui n'est pas absolument semi-simple, et un \'el\'ement semi-simple et non absolument semi-simple d'une alg\`ebre $A$ sur un corps (non parfait) $k$ ne peut admettre de d\'ecomposition de Jordan sur $A.$

Notons qu'un endomorphisme $u$ sur un $k$-espace vectoriel $E$ est semi-simple si et seulement si tout sous-espace vectoriel $F$ de $E$ invariant pour $u$ admet un suppl\'ementaire invariant pour $u,$ et que $u$ est absolument semi-simple si et seulement si il existe une extension alg\'ebrique $\tilde k$ de $k$ telle que l'extension de $u$ au $\tilde k$-espace vectoriel $\tilde E$ associ\'e \`a $E$ soit diagonalisable, voir \cite{Bou}, Chapitre VII.

Si $A$ est de dimension finie, le th\'eor\`eme de Cayley-Hamilton montre
que le polynôme caract\'eristique de $u$ annule ce dernier et on peut
le calculer explicitement. On notera que si $k$ est de caract\'eristique
nulle alors, les notations \'etant celles du th\'eor\`eme 1, on a $\tilde{p}=\frac{p}{pgcd(p,\, p')},$
et on peut calculer $pgcd(p,\, p')$ par l'algorithme d'Euclide \'etendu.

Dans le cas des alg\`ebres de matrices, on utilisera les termes \textquotedbl{}partie
diagonalisable\textquotedbl{} (sur $\tilde{k}$) et \textquotedbl{}partie
nilpotente\textquotedbl{} pour d\'esigner le terme absolument semi-simple et le
terme nilpotent de la d\'ecomposition de Jordan-Chevalley d'une matrice
$U\in{\mathcal{M}}_{n}(k).$

\section{Un exemple avec Maple}

On applique avec l'aide de Maple la m\'ethode du paragraphe pr\'ec\'edent \`a la matrice $U\in\mathcal{M}_{15}(\mathbb{C})$ d\'efinie ci-dessous 
\footnote{
Sans livrer totalement le secret de fabrication de la matrice donn\'ee ici, nous donnons un moyen simple de fabriquer une famille de matrices 15x15 non diagonalisables dont les valeurs propres ne sont pas calculables. On part d'une matrice $5\times 5$, not\'ee $A,$ \`a valeurs propres distinctes et dont les valeurs propres ne sont pas calculables.

Cela doit marcher si on prend $A$ au hasard. On note respectivement $0_5$ et $I_5$ la matrice nulle et la matrice unit\'e \`a 5 lignes et 5 colonnes, et on consid\`ere une matrice inversible $15\times 15$ not\'ee $P$.
La matrice
$U=P^{-1}\left[
\begin{matrix}
  A & I_5 & 0_5 \\
  0_5 & A & I_5 \\
  0_5 & 0_5 & A
 \end{matrix}\right]P$
est bien une matrice $15\times 15$

non diagonalisable dont les valeurs propres ne sont pas calculables. Nous laissons au lecteur le soin de trouver comment obtenir la matrice pr\'esent\'ee ci-dessus par ce proc\'ed\'e.}, dont les valeurs propres ne sont pas calculables.

\begin{scriptsize} \[U=\] \[
\left(\begin{array}{ccccccccccccccc}
-239 & -219 & -201 & -182 & -164 & -149 & -135 & -120 & -105 & -90 & -75 & -60 & -45 & -30 & -15\\
22 & 21 & 21 & 17 & 11 & 10 & 10 & 8 & 7 & 6 & 5 & 4 & 3 & 2 & 1\\
612 & 560 & 518 & 478 & 440 & 400 & 360 & 321 & 280 & 240 & 200 & 160 & 120 & 80 & 40\\
-416 & -392 & -379 & -356 & -330 & -300 & -270 & -240 & -209 & -180 & -150 & -120 & -90 & -60 & -30\\
183 & 194 & 201 & 190 & 176 & 150 & 125 & 102 & 82 & 67 & 55 & 44 & 33 & 22 & 11\\
-326 & -328 & -324 & -311 & -297 & -275 & -234 & -197 & -162 & -132 & -109 & -88 & -66 & -44 & -22\\
17 & 17 & 17 & 17 & 17 & 17 & 16 & 16 & 12 & 6 & 5 & 5 & 3 & 2 & 1\\
412 & 412 & 412 & 412 & 412 & 412 & 360 & 318 & 278 & 240 & 200 & 160 & 121 & 80 & 40\\
-266 & -266 & -266 & -266 & -266 & -266 & -242 & -229 & -206 & -180 & -150 & -120 & -90 & -59 & -30\\
128 & 128 & 128 & 128 & 128 & 128 & 139 & 146 & 135 & 121 & 95 & 70 & 47 & 27 & 12\\
-217 & -217 & -217 & -217 & -217 & -217 & -219 & -215 & -202 & -188 & -166 & -125 & -88 & -53 & -23\\
12 & 12 & 12 & 12 & 12 & 12 & 12 & 12 & 12 & 12 & 12 & 11 & 11 & 7 & 1\\
212 & 212 & 212 & 212 & 212 & 212 & 212 & 212 & 212 & 212 & 212 & 160 & 118 & 78 & 40\\
-116 & -116 & -116 & -116 & -116 & -116 & -116 & -116 & -116 & -116 & -116 & -92 & -79 & -56 & -30\\
33 & 33 & 33 & 33 & 33 & 33 & 33 & 33 & 33 & 33 & 33 & 44 & 51 & 40 & 26\end{array}\right).\]
 \end{scriptsize}
 
 Son polynome caract\'eristique $p_{U}$ , calcul\'e
sous Maple, est :\[
p=p_{U}=\left(x^{5}-9x^{4}-245x^{3}-1873x^{2}-5634x+43486\right)^{3}.\]
 Il est possible de calculer, avec la commande ``fsolve'', des
valeurs approch\'ees de chacune des valeurs propres complexes de $U,$
c'est-\`a-dire des racines de $p$. On trouve 5 valeurs propres $(\lambda_{i})_{i=1}^{5}$
, chacune ayant un ordre de muliplicit\'e 3 dans $p$ : \[
\left\{ \begin{array}{c}
\lambda_{1}\approx-9.105387869\\
\lambda_{2}\approx-4.140449458-6.991743391.i\\
\lambda_{3}\approx-4.140449458+6.991743391.i\\
\lambda_{4}\approx3.107109622\\
\lambda_{5}\approx\,23.27917716\end{array}\right..\]
 Toujours avec Maple, on calcule d'abord le pgcd $h$ de $p$ et de
son polynome d\'eriv\'e , puis on divise $p$ par $h$ pour trouver $\widetilde{p}$.
On trouve successivement : {\small \[
h=pgcd(p,\, p'
)=x^{10}-18x^{9}-409x^{8}+664x^{7}+82471x^{6}+1106154x^{5}+5486041x^{4}+\]
 \[
-203176x^{3}-131156600x^{2}-490000248x+1891032196,\]
 et\[
\widetilde{p}=\frac{p}{h}=x^{5}-9x^{4}-245x^{3}-1873x^{2}-5634x+43486.\]
 } En particulier, on se rend compte que $p=\widetilde{p}^{3},$ ce
qui signifie que deux it\'erations dans l'algorithme de Newton seront
suffisantes pour trouver la partie diagonalisable. Pour pouvoir faire
tourner l'algorithme, on a besoin du polynôme d\'eriv\'e de $\widetilde{p}$
qui se calcule sous Maple \`a l'aide de la commande {}``diff''. On
trouve\[
\tilde{p}'=5x^{4}-36x^{3}-735x^{2}-3746x-5634.\]
 Puis, on en d\'eduit $q$ tel que $q\widetilde{p}'\equiv1\;(mod\,\overline{p}),$
\`a l'aide de la commande {}``gcdex''. Enfin on calcule $q(U),$
$\widetilde{p}(U)$ et \[
D_{1}=U-\widetilde{p}(U)q(U).\]
 De même, on calcule $D_{2}$ \`a l'aide de la formule\[
D_{2}=D_{1}-\widetilde{p}(D_{1})q(D_{1}).\]
 On trouve ainsi la partie diagonalisable $D$ de la matrice $U.$

\begin{scriptsize}\[
D_{2}=D=\]
 \[
\left(\begin{array}{ccccccccccccccc}
-240 & -220 & -202 & -183 & -165 & -150 & -135 & -120 & -105 & -90 & -75 & -60 & -45 & -30 & -15\\
22 & 21 & 21 & 17 & 11 & 10 & 9 & 8 & 7 & 6 & 5 & 4 & 3 & 2 & 1\\
612 & 560 & 518 & 478 & 440 & 400 & 360 & 320 & 280 & 240 & 200 & 160 & 120 & 80 & 40\\
-416 & -392 & -379 & -356 & -330 & -300 & -270 & -240 & -210 & -180 & -150 & -120 & -90 & -60 & -30\\
183 & 194 & 201 & 190 & 176 & 150 & 125 & 102 & 82 & 66 & 55 & 44 & 33 & 22 & 11\\
-326 & -328 & -324 & -311 & -297 & -275 & -234 & -197 & -162 & -132 & -110 & -88 & -66 & -44 & -22\\
17 & 17 & 17 & 17 & 17 & 17 & 16 & 16 & 12 & 6 & 5 & 4 & 3 & 2 & 1\\
412 & 412 & 412 & 412 & 412 & 412 & 360 & 318 & 278 & 240 & 200 & 160 & 120 & 80 & 40\\
-266 & -266 & -266 & -266 & -266 & -266 & -242 & -229 & -206 & -180 & -150 & -120 & -90 & -60 & -30\\
128 & 128 & 128 & 128 & 128 & 128 & 139 & 146 & 135 & 121 & 95 & 70 & 47 & 27 & 11\\
-216 & -216 & -216 & -216 & -216 & -216 & -218 & -214 & -201 & -187 & -165 & -124 & -87 & -52 & -22\\
12 & 12 & 12 & 12 & 12 & 12 & 12 & 12 & 12 & 12 & 12 & 11 & 11 & 7 & 1\\
212 & 212 & 212 & 212 & 212 & 212 & 212 & 212 & 212 & 212 & 212 & 160 & 118 & 78 & 40\\
-116 & -116 & -116 & -116 & -116 & -116 & -116 & -116 & -116 & -116 & -116 & -92 & -79 & -56 & -30\\
33 & 33 & 33 & 33 & 33 & 33 & 33 & 33 & 33 & 33 & 33 & 44 & 51 & 40 & 26\end{array}\right),\]
 \end{scriptsize}

puis, par diff\'erence, sa partie nilpotente $N$:

\begin{scriptsize} \[
N=A-D=\left(\begin{array}{ccccccccccccccc}
1 & 1 & 1 & 1 & 1 & 1 & 0 & 0 & 0 & 0 & 0 & 0 & 0 & 0 & 0\\
0 & 0 & 0 & 0 & 0 & 0 & 1 & 0 & 0 & 0 & 0 & 0 & 0 & 0 & 0\\
0 & 0 & 0 & 0 & 0 & 0 & 0 & 1 & 0 & 0 & 0 & 0 & 0 & 0 & 0\\
0 & 0 & 0 & 0 & 0 & 0 & 0 & 0 & 1 & 0 & 0 & 0 & 0 & 0 & 0\\
0 & 0 & 0 & 0 & 0 & 0 & 0 & 0 & 0 & 1 & 0 & 0 & 0 & 0 & 0\\
0 & 0 & 0 & 0 & 0 & 0 & 0 & 0 & 0 & 0 & 1 & 0 & 0 & 0 & 0\\
0 & 0 & 0 & 0 & 0 & 0 & 0 & 0 & 0 & 0 & 0 & 1 & 0 & 0 & 0\\
0 & 0 & 0 & 0 & 0 & 0 & 0 & 0 & 0 & 0 & 0 & 0 & 1 & 0 & 0\\
0 & 0 & 0 & 0 & 0 & 0 & 0 & 0 & 0 & 0 & 0 & 0 & 0 & 1 & 0\\
0 & 0 & 0 & 0 & 0 & 0 & 0 & 0 & 0 & 0 & 0 & 0 & 0 & 0 & 1\\
-1 & -1 & -1 & -1 & -1 & -1 & -1 & -1 & -1 & -1 & -1 & -1 & -1 & -1 & -1\\
0 & 0 & 0 & 0 & 0 & 0 & 0 & 0 & 0 & 0 & 0 & 0 & 0 & 0 & 0\\
0 & 0 & 0 & 0 & 0 & 0 & 0 & 0 & 0 & 0 & 0 & 0 & 0 & 0 & 0\\
0 & 0 & 0 & 0 & 0 & 0 & 0 & 0 & 0 & 0 & 0 & 0 & 0 & 0 & 0\\
0 & 0 & 0 & 0 & 0 & 0 & 0 & 0 & 0 & 0 & 0 & 0 & 0 & 0 & 0\end{array}\right).\]
\end{scriptsize}

\section{D\'ecomposition de Jordan via un syst\`eme de congruences}

On reprend ici les notations du th\'eor\`eme 1, et on d\'esigne par $p=p_{1}^{m_{1}}\dots p_{r}^{m_{r}}$
un polynôme annulateur d'un \'el\'ement s\'eparable $u$ d'une $k$-alg\`ebre
$A$ dont les facteurs irr\'eductibles $p_{1},\dots,\, p_{r}$ sont
\`a racines simples dans $\tilde{k}.$ Notons $\lambda_{1},...,\lambda_{s}$
les racines distinctes de $\tilde{p}:=p_{1}\dots p_{r},$ et $n_{1},\dots,\, n_{s}$
leurs ordres de multiplicit\'e (notons que pour tout $1\leq i\leq r,$
toutes les racines de $p_{i}$ sont en fait d'ordre $m_{i}$). On
consid\`ere le syst\`eme d'\'equations de congruence

\[
\left(\mathcal{S}\right):\:\left\{ \begin{array}{c}
h\equiv\lambda_{1}\:\:{\rm mod}\ \left(x-\lambda_{1}\right)^{n_{1}}\\
h\equiv\lambda_{2}\:\:{\rm mod}\ \left(x-\lambda_{2}\right)^{n_{2}}\\
\vdots\\
\vdots\\
h\equiv\lambda_{s}\:\:{\rm mod}\ \left(x-\lambda_{s}\right)^{n_{s}}\end{array}\right..\]
 Il r\'esulte du th\'eor\`eme chinois que ce syst\`eme poss\`ede une solution
(unique modulo $p$) dans $\tilde{k}[x].$ Le fait que ce syst\`eme
poss\`ede une solution dans $k[x]$ est un peu moins \'evident. On pose
\[
m:=\max_{1\le j\le s}n_{j}=\max_{1\le j\le r}m_{j},\]
 et on consid\`ere le syst\`eme un peu plus compliqu\'e

\[
\left(\mathcal{S}'\right):\:\left\{ \begin{array}{c}
h\equiv\lambda_{1}\:\:{\rm mod}\ \left(x-\lambda_{1}\right)^{m}\\
h\equiv\lambda_{2}\:\:{\rm mod}\ \left(x-\lambda_{2}\right)^{m}\\
\vdots\\
\vdots\\
h\equiv\lambda_{s}\:\:{\rm mod}\ \left(x-\lambda_{s}\right)^{m}\end{array}\right..\]
 Le syst\`eme $S'$ admet une unique solution modulo $\tilde{p}^{m}.$
De plus, comme $\tilde{p}\in k[x]$ est \`a racines simples, il
existe une extension galoisienne $k_{1}\subset\tilde{k}$ de $k$
contenant $\lambda_{1},...,\lambda_{s}.$ Ceci signifie que l'ensemble
des \'el\'ements de $k_{1}$ invariants pour tous les automorphismes $\sigma\in G$
est r\'eduit \`a $k,$ où $G:=Gal(k_{1}/k)$ d\'esigne le groupe de Galois
form\'e des automorphismes de $k_{1}$ laissant fixes tous les \'el\'ements
de $k.$ Soit $h_{0}\in\tilde{k}(x)$ l'unique solution de degr\'e strictement
inf\'erieur \`a $sm$ de ${\mathcal{S}}',$ soit $\sigma\in G,$ et soit
$\sigma(h_{0})$ le polynôme obtenu en appliquant $\sigma$ \`a tous
les coefficients de $h_{0}.$ Comme $\sigma(\lambda_{j})\in\{\lambda_{1},\dots,\,\lambda_{s}\},$
et comme $\sigma(h_{0})\equiv\sigma(\lambda_{j})\:\:{\rm mod}\ \left(x-\sigma(\lambda_{j})\right)^{n_{j}}$
pour $1\le j\le s,$ on voit que $\sigma(h_{0})$ est aussi une solution
de ${\mathcal{S}}'$ de degr\'e strictement inf\'erieur \`a $sm.$ Donc
$\sigma(h_{0})=h_{0}$ pour tout $\sigma\in G,$ et $h_{0}\in k[x].$
Notons qu'en remplaçant $h_{0}$ par le reste de sa division par $p,$
on obtient la solution $h\in k[x]$ du syst\`eme $\mathcal{S}$ de degr\'e
strictement inf\'erieur \`a celui de $p.$

On a $\tilde{p}(h(u))=(h(u)-\lambda_{1})\dots(h(u)-\lambda_{s})\in p(u)B=\{0\},$
donc $\tilde{p}(d)=0$ et $h(u)$ est absolument semi-simple. De plus, $x-h$
est divisible par $x-\lambda_{j}$ pour $1\le j\le n,$ donc $x-h$
est divisible par $\tilde{p}=ppcm(x-\lambda_{1},\dots,\, x-\lambda_{s}),$
et $(x-h)^{m}$ est divisible par $\tilde{p}^{m}=p.$ Donc $u-h(u)$
est nilpotent, et $(h(u),\, u-h(u))$ est la d\'ecomposition de Jordan-Chevalley
de $u.$

\textit{Algorithme de Newton.} En fait on peut calculer $h$ en utilisant
l'algorithme de Newton. Il s'agit d'une simple application du th\'eor\`eme
1. En effet posons ${\mathcal{A}}=k[x]/pk[x],$ soit $\pi=k[x]\to{\mathcal{A}}$
la surjection canonique et soit $v=\pi(x).$ On a $p(v)=0,$ et $v$
est un \'el\'ement s\'eparable de ${\mathcal{A}}.$ Posons $\alpha=d^{o}(p).$
En identifiant ${\mathcal{A}}$ \`a l'ensemble $k_{\alpha-1}[x]$ des
polynômes de degr\'e inf\'erieur ou \'egal \`a $\alpha-1,$ l'algorithme de
Newton donne la suite $(h_{n})_{n\ge0}$ de polynômes telle que $h_{0}=x$
et telle que $h_{n+1}$ est le reste de la division euclidienne de
$h_{n}-\tilde{p}(h_{n})q(h_{n})$ par $p.$ On a $h_{n}(\lambda_{j})=\lambda_{j}$
pour $1\le j\le s,$ $n\ge0,$ et $\tilde{p}(h_{N})=0,$ ce qui signifie
que $\prod_{1\le j\le s}(h_{N}-\lambda_{j})$ est divisible par $p.$

On a $h_{N}(\lambda_{j})=\lambda_{j}$ pour $1\le j\le s,$ $n\ge0,$
et $x-\lambda_{j}$ est premier avec $\prod_{l\neq j}(h_{N}-\lambda_{l}).$
Donc $(x-\lambda_{j})^{n_{j}}$ divise $h_{N}-\lambda_{j}$ pour $1\le j\le s,$
et $h_{N}$ est solution du syst\`eme ${\mathcal{S}}.$

Consid\'erons de nouveau l'exemple de la matrice $U\in\mathcal{M}_{15}(\mathbb{C})$
\'etudi\'ee au cours de la partie III. Elle donne $\alpha=15,$ $n_{1}=n_{2}=n_{3}=3,$
et $m=2.$ Voici le polynôme $h_{2}$ de degr\'e 14 trouv\'e en utilisant
Maple:

\begin{scriptsize} \[
-\frac{164777455994373388396328621588559}{6153698372078022255427531347585869793027770912}x^{14}+\frac{4696495785239673852290048403064815}{6153698372078022255427531347585869793027770912}x^{13}\]

\smallskip{}

\[
+\frac{1904832784435747945567751656823474}{192303074127438195482110354612058431032117841}x^{12}-\frac{93156612626828282836743908787001857}{769212296509752781928441418448233724128471364}x^{11}\]

\smallskip{}
 \[
-\frac{947888046084076156531978456271155215}{192303074127438195482110354612058431032117841}x^{10}-\frac{134225331044956775638702555509762875765}{3076849186039011127713765673792934896513885456}x^{9}+\]

\smallskip{}
 \[
+\frac{1437307520588625416353772759457723214755}{3076849186039011127713765673792934896513885456}x^{8}+\frac{5160755602829187090655465367051014084185}{384606148254876390964220709224116862064235682}x^{7}+\]

\smallskip{}
 \[
+\frac{883915408036735072058939076768746842566791}{6153698372078022255427531347585869793027770912}x^{6}+\frac{3913459854154300022640900141580672168912059}{6153698372078022255427531347585869793027770912}x^{5}+\]

\smallskip{}

\[
-\frac{4125757581287724475079189681241424637007729}{3076849186039011127713765673792934896513885456}x^{4}-\frac{1202955633870054250571522565744213358793813}{41579043054581231455591428024228849952890344}x^{3}+\]

\smallskip{}
 \[
-\frac{201785886295705753509895176779863737263724099}{1538424593019505563856882836896467448256942728}x^{2}+\frac{1811980652583749695251762135187625041505888535}{1538424593019505563856882836896467448256942728}x+\]

\smallskip{}
 \[
+\frac{1040926769591787693101439601278755401419987857}{769212296509752781928441418448233724128471364}.\]

\end{scriptsize}

\bigskip{}

Une valeur approch\'ee de $h_{2}$ se calcule \`a l'aide de la commande
{}``evalf'' et vaut :

\smallskip{}
 \begin{scriptsize} \[
-0.2677697964.10^{-13}x^{14}+0.7631988930.10^{-12}x^{13}+0.9905368352.10^{-11}x^{12}-0.1211065047.10^{-9}x^{11}\]
 \[
-0.4929136211.10^{-8}x^{10}-0.4362428021.10^{-7}x^{9}+0.4671361623.10^{-6}x^{8}+0.00001341828680.x^{7}+0.0001436397032.x^{6}\]
 \[
+0.0006359524984.x^{5}-0.001340903415.x^{4}-.02893177778.x^{3}-0.1311639759.x^{2}+1.177815709.x+1.353237298\,.\]

\end{scriptsize} 

\section{La terminologie \textquotedbl{}d\'ecomposition de Jordan-Chevalley\textquotedbl{}}

La d\'efinition 1 de notre deuxi\`eme partie donne le nom de \textquotedbl{}d\'ecomposition
de Jordan-Chevalley\textquotedbl{} au couple $(d,\, s)$ d'\'el\'ements
de $A$ tels que $u=d+s$.

Le nom de Jordan, tr\`es classique dans le contexte de d\'ecomposition
d'un endomorphisme, ne peut \'etonner ici.

Certes, la d\'ecomposition d\'ecrite dans la d\'efinition ci-dessus est
fort \'eloign\'ee du th\'eor\`eme \'enonc\'e par Camille Jordan dans le {\em
Trait\'e des substitutions et des \'equations alg\'ebriques} \cite{Jor70}
de 1870. A cette \'epoque, le th\'eor\`eme dû \`a C.~Jordan prend le nom
de <<forme canonique d'une substitution lin\'eaire>> et est \'ecrit
en termes de substitutions lin\'eaires de la forme $|x,x'\cdots ax+bx'+\cdots,a'x+b'x'+\cdots,\cdots|$,
l'\'ecriture matricielle n'apparaissant jamais. La d\'emonstration donn\'ee
par C.~Jordan est, elle aussi, bien loin de nos m\'ethodes actuelles.
On peut la trouver \'etudi\'ee et comment\'ee \`a l'aide d'exemples par F.~Brechenmacher (\cite{BreT06}, p.182, encart
5).

C'est seulement \`a partir des ann\'ees 1930-1940 que les versions successives
de ce th\'eor\`eme sont d\'efinitivement rattach\'ees au nom du math\'ematicien
Camille Jordan. On peut citer par exemple en 1932 H.W.~Turnbull et
A.C.~Aitken page 114:
\begin{quotation}
<<The Classical Form C is first found in C.~Jordan, {\em Trait\'e
des substitutions}>>\cite{Tur32}
\end{quotation}
ou en 1943 C.C.~Mac Duffee th\'eor\`eme 65:
\begin{quotation}
<<This is the familiar Jordan normal form of a matrix with complex
elements.>> \cite{Duf43}.
\end{quotation}
Dans ce cadre, il est toujours question d'un th\'eor\`eme de r\'eduction
(et non de d\'ecomposition), tel qu'on le trouve encore dans les ann\'ees
1970, par exemple dans l'ouvrage de Serge Lang (page 398, {\em Algebra}
\cite{Lan65}).

Pourtant, en parall\`ele, une autre histoire du th\'eor\`eme de Jordan se
dessine vers 1950. C'est l\`a que se trouve la source de la d\'enomination
\textquotedbl{}d\'ecomposition de Jordan-Chevalley\textquotedbl{}.

Claude Chevalley, qui fait partie du groupe Bourbaki depuis sa cr\'eation
en 1935, oriente ses recherches en direction des groupes de Lie et
des groupes alg\'ebriques et publie depuis les États-Unis \cite{Die99}
deux livres sur ces sujets :

{\em Theory of Lie groups I} en 1946 et {\em Th\'eorie des groupes
de Lie Tome II} (sous titre : groupes alg\'ebriques) en 1951.

Claude Chevalley n'est \'evidemment pas seul \`a travailler sur ces questions.
On peut trouver une \'etude approfondie de ce domaine math\'ematique dans
le livre \emph{Essays in the History of Lie Groups and Algebraic Groups}
qu'Armand Borel a consacr\'e au <<first century of the history of Lie
groups and algebraic groups>> (\cite{Bor01} page IX). Compl\`etement
contemporains des travaux de Claude Chevalley, citons ceux d' E.R.~Kolchin,
qui publie en 1946 sur les groupes alg\'ebriques, et fait appel dans
ses travaux \`a la <<forme normale de Jordan>> (\cite{Kol48} p. 9,
\cite{Kol48-2} p. 776).

Revenons au deuxi\`eme tome de l'ouvrage de Claude Chevalley \cite{Che51};
au paragraphe {\em Espaces vectoriels \`a op\'erateurs}, on trouve,
page 71 le th\'eor\`eme 7 :
\begin{quotation}
<<Soit $X$ un endomorphisme d'un espace vectoriel $V$ de dimension
finie sur un corps parfait $k$. Il est alors possible, d'une mani\`ere
et d'une seule, de repr\'esenter $X$ comme somme d'un endomorphisme
semi-simple $S$ et d'un endomorphisme nilpotent $N$ qui commutent
entre eux; $S$ et $N$ peuvent être repr\'esent\'es comme polynômes en
$X$ \`a coefficients dans $k$.>>
\end{quotation}
Clairement, ici, Claude Chevalley ouvre la porte de la \textquotedbl{}d\'ecomposition\textquotedbl{}
de Jordan.

Dans ce livre, le but de Claude Chevalley n'est pas d'\'ecrire la matrice
sous une <<forme aussi simple que possible>> comme le recherchait
Camille Jordan (\cite{Jor70}, p.114). Il utilise ce th\'eor\`eme de d\'ecomposition
afin de pouvoir d\'emontrer le th\'eor\`eme final de l'ouvrage, th\'eor\`eme
18 page 184, sur les groupes alg\'ebriques lin\'eaires:
\begin{quotation}
<<Tout groupe alg\'ebrique d'automorphismes de $V$ qui contient $s$
contient aussi $u$ et $v$>>.
\end{quotation}
(avec $s=uv$, $u$ et $v$ \'etant les composantes semi-simple et unipotente
de l'automorphisme $s$, $V$ \'etant un espace vectoriel de dimension
finie sur un corps K de caract\'eristique $0$).

Donnons maintenant quelques rep\`eres pour situer le devenir de ces
deux th\'eor\`emes.

D\`es 1956, Armand Borel \'ecrit un important article sur les groupes
alg\'ebriques \cite{Bor56}. A la premi\`ere page de l'introduction, on
lit :
\begin{quotation}
<<On sait que tout automorphisme $g$ d'un espace vectoriel $V$
de dimension finie sur un corps parfait $k$ s'\'ecrit d'une façon et
d'une seule comme produit $g_{s}.g_{u}$ d'un automorphisme $g_{s}$
semi-simple (i.e. \`a diviseurs \'el\'ementaires) et d'un automorphisme
$g_{u}$ unipotent (i.e. \`a valeurs propres \'egales \`a 1) commutant entre
eux.>>
\end{quotation}
Comme chez Claude Chevalley, dans ce contexte des groupes alg\'ebriques,
la \textquotedbl{}d\'ecomposition de Jordan\textquotedbl{} (\'ecrite ici
sous forme multiplicative et pour les automorphismes) prime. Un des
th\'eor\`emes centraux de l'article d'Armand Borel (th\'eor\`eme 8-4 page
47) n'est autre que le th\'eor\`eme 18 de Claude Chevalley sur les groupes
alg\'ebriques donn\'e plus haut. Et Armand Borel cite en r\'ef\'erence Claude
Chevalley \cite{Che51}.

Un peu plus tard, le th\'eor\`eme de \textquotedbl{}d\'ecomposition de Jordan\textquotedbl{}
apparait \`a nouveau dans le texte du S\'eminaire dirig\'e par Claude Chevalley
\`a l'École Normale Sup\'erieure (ann\'ees universitaires 56/57 et 57/58
dans une partie r\'edig\'ee par Alexandre Grothendieck). On peut y lire
pages 47-48 de \cite{Che05}:  \begin{quotation}<<Rappelons le fait bien connu:
$x$ peut se mettre de façon unique sous la forme $x_{s}$+$x_{n}$,
somme d'un endomorphisme semi-simple $x_{s}$ et d'un endomorphisme
nilpotent $x_{n}$ qui commutent (appel\'es partie semi-simple et partie
nilpotente de $x$).>>
\end{quotation}
C'est bien le même \'enonc\'e que celui de Claude Chevalley dans \cite{Che51},
$x$ \'etant un endomorphisme de $V$, espace vectoriel de dimension
finie. La forme multiplicative pour un automorphisme apparaît quelques
lignes plus loin. Le texte du S\'eminaire Chevalley, suscitant un grand
int\'erêt dans la communaut\'e math\'ematique ({\em Avertissement au lecteur}
\cite{Che05}), a \'et\'e r\'e\'edit\'e en 2005.

Ainsi, peu \`a peu, c'est sous cette forme que le \textquotedbl{}th\'eor\`eme
de Jordan\textquotedbl{} va devenir r\'ef\'erence et prendre le nom de
Jordan-Chevalley quand il en est question dans la th\'eorie des groupes
alg\'ebriques. C'est ce que l'on peut voir par exemple dans les ouvrages
de James E.~Humphreys (\cite{Hum72}, \cite{Hum}). En particulier,
dans le deuxi\`eme, {\em Linear algebraic groups}, c'est tout le
paragraphe reprenant les th\'eor\`emes 7 et 18 de Claude Chevalley qui
prend le nom de << Jordan-Chevalley Decomposition>>.

Une autre raison de notre d\'enomination de \textquotedbl{}th\'eor\`eme
de Jordan-Chevalley\textquotedbl{} se trouve dans le livre {\em
Th\'eorie des groupes de Lie Tome II} de Claude Chevalley. Revenons
sur quelques \'el\'ements de la d\'emonstration du th\'eor\`eme 7 telle qu'elle
y est pr\'esent\'ee. Pour la lire, nous avons besoin de connaître la construction
du polynôme $f$, donn\'ee par Claude Chevalley dans la proposition
5 qui pr\'ec\`ede le th\'eor\`eme 7 :
\begin{quotation}
<<Soit $X$ un endomorphisme d'un espace vectoriel $V$ de dimension
finie sur un corps parfait $k$. Pour que $X$ soit semi-simple, il
faut et il suffit qu'il existe un polynôme $f$ \`a coefficients dans
$k$, relativement premier \`a son polynôme d\'eriv\'e, tel que $f(X)=0$.>>
\end{quotation}
Un endomorphisme semi-simple est donc caract\'eris\'e comme dans notre
d\'efinition 1, Claude Chevalley se plaçant, comme nous l'avons vu,
dans le cas où $X$ est un endomorphisme d'un espace vectoriel $V$
de dimension finie sur un corps parfait $k$, cadre moins g\'en\'eral
que celui de notre deuxi\`eme partie.

En suivant la d\'emonstration propos\'ee par Claude Chevalley dans \cite{Che51},
on y retrouve le fil conducteur des d\'emonstrations ult\'erieures, comme
celle de notre th\'eor\`eme 1.

Dans le texte de C. Chevalley, on suppose connu un polynôme $F$ \`a
coefficients dans $k$, tel que $F(X)=0$ (par exemple le polynôme
caract\'eristique de $X$). Et on lit :
\begin{quotation}
<<Écrivons $F=c{F_{1}}^{e_{1}}\dots{F_{h}}^{e_{h}}$ où $c\in K$,
$F_{1},...F_{h}$ sont des polynômes irr\'eductibles relativement premiers
entre eux deux \`a deux, et $e_{1},\dots,\, e_{h}$ des exposants $>0$.
Soit $f=F_{1}\dots F_{h}$ et soit $e$ le plus grand des exposants
$e_{i}$ ($1\leq i\leq h$). Le polynôme $f^{e}$ est alors divisible
par $F$ d'où $(f(X))^{e}=0$.>>
\end{quotation}
Le polynôme $f$ est tout simplement notre polynôme $\tilde{p}$.

Puis, dans cette d\'emonstration, Claude Chevalley construit explicitement
le polynôme \`a une variable $U$, not\'e $s^{m}(U)$ dans son texte,
tel que $s^{m}(U)(X)$ soit la partie diagonalisable (ou semi-simple)
de l'endomorphisme $X$. C'est ce sch\'ema de d\'emonstration que l'on
retrouve dans les constructions de polynômes amenant \`a la d\'ecomposition
de Jordan, comme celle expos\'ee dans le th\'eor\`eme 1 de l'article. Ces
constructions s'appuient sur la \textquotedbl{}m\'ethode de Newton\textquotedbl{}
et aboutissent pour nous \`a la construction de $d_{N}$, partie semi-simple
de la d\'ecomposition de $u$, $u$ \'etant un \'el\'ement s\'eparable de l'alg\`ebre
$A$.

Ce rapide d\'etour par la d\'emonstration de 1951 met en lumi\`ere le rôle
de Claude Chevalley et permet de confirmer la d\'enomination \textquotedbl{}th\'eor\`eme
de Jordan-Chevalley\textquotedbl{}.

Certes, en 1951, dans son livre, Claude Chevalley ne s'arrête pas
sur le fait que $f$, que nous avons not\'e $\tilde{p}$, peut se calculer
sans connaître les valeurs propres de l'endomorphisme $X$. La m\'ethode,
que nous avons d\'etaill\'ee pr\'ec\'edemment, s'appuie sur le calcul de $f$
(ou $\tilde{p}$) avec $\tilde{p}=\frac{p}{pgcd(p,p')}$ où $p$ est
le polynôme caract\'eristique de $X$. Si la forme du polynôme $f$
se trouve bien chez Claude Chevalley, il ne l'utilise pas dans ce
but.

On peut noter par ailleurs que la d\'emonstration qui nous int\'eresse
aujourd'hui n'est reprise, ni par Armand Borel, ni par Alexandre Grothendieck,
ni par J.E.~Humphreys quand ils travaillent sur les groupes alg\'ebriques. L'ouvrage {\em Th\'eorie des groupes de Lie Tome II} n'\'etait pas
destin\'e \`a l'enseignement en premier cycle.

\section{Quelques aspects de la d\'ecomposition multiplicative de Jordan-Chevalley dans les groupes de Lie semi-simples}
La d\'ecomposition de Jordan joue bien s\^ur un r\^ole important dans la th\'eorie des groupes et alg\`ebres de Lie, pour laquelle nous renvoyons \`a \cite{He}, \cite{Va}. Soit $\mathfrak {g}$ une alg\`ebre de Lie semi-simple r\'eelle ou complexe, soit $Der(\mathfrak {g})$ l'ensemble des d\'erivations de $\mathfrak {g},$ c'est-\`a-dire l'ensemble des endomorphismes $D$ de $\mathfrak {g}$ v\'erifiant l'identit\'e $D([Y,Z])=[D(Y),Z] +[Y,D(Z)],$ et pour $X \in \mathfrak{g}$ soit $ad(X): Y \to [X,Y]$ la d\'erivation associ\'ee \`a $X.$ Comme
$\mathfrak {g}$ est semi-simple, la repr\'esentation adjointe $X \to ad(X)$ est un isomorphisme de $\mathfrak {g}$ sur $Der(\mathfrak {g}).$  De plus, Claude Chevalley a montr\'e dans \cite{Che51}, corollaire du th\'eor\`eme 16 que la
 partie  semi-simple et la partie nilpotente de la d\'ecomposition de Jordan d'une d\'erivation sur une alg\`ebre r\'eelle ou complexe $A$, non n\'ecessairement associative, sont des d\'erivations sur $A.$ On  d\'eduit alors du th\'eor\`eme de d\'ecomposition de Jordan que pour tout $X\in \mathfrak {g}$ il existe un unique couple $(S,N)$ d'\'el\'ements de $\mathfrak {g}$ tels que $ad(S)$ soit semi-simple, $ad(N)$ nilpotent et $[S,N]=0,$ voir  \cite{Va}, th\'eor\`eme 3.10.6. 

Soit $G$ un groupe de Lie semi-simple r\'eel ou complexe d'unit\'e $1_G$, soit $G_0$ la composante connexe de $1_G$ dans $G,$ soit ${\mathfrak g}$ l'alg\`ebre de Lie de $G,$ identifi\'ee \`a l'espace tangent \`a $G$ en $1_G,$ et soit  $\mathfrak {Aut(g)}$ l'ensemble des isomorphismes lin\'eaires $\theta$ de $\mathfrak {g}$ v\'erifiant  $\theta([X,Y])=[\theta(X),\theta(Y)]$ pour
$X,Y\in \mathfrak{g}.$  Pour $g \in G$ soit $Ad(g)\in \mathfrak{Aut(g)}$ la diff\'erentielle en $1_G$ de  l'automorphisme de conjugaison $\phi_g: a \to gag^{-1}.$ La repr\'esentation adjointe $g \to Ad(g)$ de $G$ n'est en g\'en\'eral pas injective, puisque son noyau est le centralisateur $\{ g \in G \ | \ ag=ga \ \forall a \in G_0\}$ de $G_0$  dans $G.$ On dira que $g \in G$ est {\it semi-simple} si $Ad(g)$ est un endomorphisme semi-simple de $G,$ et on dira que $g \in G$ est {\it unipotent} si $g=\exp(X),$ ce qui donne $Ad(g)=\exp(ad(X)),$ o\`u $X\in \mathfrak {g}$ est tel que $ad(X)$ soit nilpotent.

Soit maintenant  $G$  un groupe de Lie semi-simple complexe. Comme $\mathfrak {Aut(g)}$ est un sous-groupe alg\'ebrique du groupe des automorphismes lin\'eaires de $\mathfrak {g},$ il r\'esulte du th\'eor\`eme 18 de \cite{Che51} cit\'e plus haut qu'il existe un couple unique
$(S, U)$ d'\'el\'ements de $\mathfrak{Aut(g)}$ v\'erifiant $Ad(g)=US=SU$ tels que $S$ (resp. $U$) soit un endomorphisme semi-simple (resp. unipotent) de $\mathfrak{g}.$ Posons $N= U-I,$   et  $N_0=\sum_{n=0}^{m-1}(-1)^{n+1}{N^n\over n},$ de sorte que $N$ et $N_0$ sont nilpotents. Comme $\exp(tN_0)$ est un polyn\^ome en $t,$ et comme $\exp(N_0)=U,$ on a $\exp(tN_0)\in \mathfrak{Aut(g)}$ pour $t >0$, ce qui implique que $N_0\in Der(\mathfrak{g}),$ voir \cite{Va}, Chap. 2, exercice 21, et il existe $X \in \mathfrak {g}$ tel que $N_0=ad(X),$ d'o\`u $U=Ad(\exp(X)),$ $S=Ad(g\exp(-X)).$ On voit donc que {\it tout \'element d'un groupe de Lie semi-simple complexe s'\'ecrit comme produit d'un \'el\'ement semi-simple et d'un \'el\'ement unipotent du groupe qui commutent entre eux} et on v\'erifie que cette d\'ecomposition est unique.

On obtient un r\'esultat analogue pour un groupe de Lie semi-simple  r\'eel en introduisant la complexifi\'ee de son alg\`ebre de Lie. Nous \'enon\c cons un r\'esultat plus pr\'ecis pour un groupe de Lie r\'eel semi-simple et connexe $G$. On dit que $g \in G$ est elliptique si les valeurs propres de $Ad(g),$ consid\'er\'e comme endomorphisme de la complexifi\'ee de l'alg\`ebre de Lie de $\mathfrak g$ de $G,$ sont de module 1, et on dira que $g \in G$ est hyperbolique s'il existe $X \in \mathfrak{g},$ dont l'adjoint $ad(X)$ est diagonalisable sur l'espace vectoriel r\'eel  $\mathfrak g$, tel que $g=\exp(X).$

On a alors le r\'esultat suivant, voir l'article de Bertram Kostant \cite{Ko}

\textbf{Th\'eor\`eme 2.}\textit{ Soit $G$ un groupe de Lie r\'eel semisimple et connexe. Alors tout \'el\'ement $g$ de $G$ se d\'ecompose de mani\`ere unique sous la forme $$g=ehu$$ o\`u $e,$ $h,$ $u$ sont des \'el\'ements respectivement elliptiques, hyperboliques et unipotents de $G$ qui commutent entre eux.}

 Autrement dit la partie semi-simple de $d$ de $g,$ dont l'adjoint est diagonalisable sur la complexifi\'ee de $\mathfrak g,$ se d\'ecompose en produit d'un \'element elliptique et d'un \'el\'ement hyperbolique de $G$ qui commutent entre eux, et on obtient pour les groupes de Lie semi-simples connexes une "d\'ecomposition de Jordan multiplicative compl\`ete" analogue \`a la d\'ecomposition de Jordan multiplicative compl\`ete usuelle dans $GL(n,\mathbb R),$ donn\'ee par exemple dans  \cite{He}, lemme 7.1. On trouve dans le th\'eor\`eme 7.2 de \cite{He}, le lien classique entre ellipticit\'e, hyperbolicit\'e et nilpotence et d\'ecompositions d'Iwasawa \footnote{La "d\'ecomposition d'Iwasawa", not\'ee $G = KAN,$ d\'ecompose un groupe de Lie semi-simple en produit d'un sous-groupe compact maximal $K,$ d'un sous groupe de Cartan $A,$ sous-groupe de Lie correspondant \`a une sous-alg\`ebre de Lie commutative maximale $\mathfrak h$, et d'un groupe nilpotent $N$; ainsi par exemple pour $GL(n, \mathbb R)$, $K$ est le groupe orthogonal, $A$ le sous-groupe form\'e des matrices diagonales et $N$ le sous-goupe des matrices triangulaires sup\'erieures dont tous les termes diagonaux sont \'egaux \`a 1.}

Pour conclure cette br\`eve discussion historique nous mentionnons un r\'esultat tr\`es r\'ecent qui concerne les groupes de Lie r\'eels semi-simples lin\'eaires, c'est-\`a-dire les groupes de Lie r\'eels semi-simples qui sont des sous-groupes ferm\'es de $GL(n,\,\mathbb R).$ Soit $G$ un groupe localement compact, soit $\Gamma$ un sous-groupe discret de $G,$ soit $\Gamma\setminus G$ l'ensemble des classes \`a droite modulo $\Gamma,$ muni de la topologie quotient, soit  ${\mathcal C}_c(G)$ l'ensemble des fonctions continues \`a support compact sur $G$, et soit $\mu$ une mesure de Haar \`a gauche sur $G.$  On peut  d\'efinir une mesure positive $G$-invariante $\nu$ sur $\Gamma\setminus G$ par la formule
$$\int_{\Gamma \setminus G}\tilde f(\eta)d\nu(\eta)=\int_{g \in G} f(g)d\mu(g) \ \ \ \ \ (f \in {\mathcal C}_c(G)),$$

o\`u $\tilde f(\Gamma x)=\sum_{\gamma \in \Gamma}f(\gamma x)$ pour $x\in G$ (l'ouvrage de Lynn Loomis \cite{Lo}, accessible en ligne, reste une bonne r\'ef\'erence pour les mesures invariantes sur les groupes quotient).

On dit que le groupe discret $\Gamma$ est un {\it r\'eseau} de G si $\Gamma$ est de covolume fini, c'est-\`a-dire si $\int_{\Gamma \setminus G}d\nu(\eta)<+\infty,$

N.T. Venkataramana a montr\'e en 2008 dans \cite{Ve} que si $\Gamma$ est un r\'eseau d'un groupe de Lie lin\'eaire semi-simple r\'eel $G$, et si $\gamma \in \Gamma,$ alors la partie semi-simple $s$ et la partie unipotente $u$ de la d\'ecomposition de Jordan multiplicative de $\gamma$ v\'erifient  $s
^m\in \Gamma$ et $u^m\in \Gamma$ pour un certain entier $m\ge 1,$ r\'esultat qui \'etait classique pour les sous-groupes "arithm\'etiques" de $G.$ Si de plus le "rang r\'eel" de $G$ est \'egal \`a 1, et si $u\neq 1,$ alors il existe $n\ge 1$
tel que $s^n=1_G$. Nous renvoyons le lecteur \`a l'ouvrage de Gregori Margulis \cite {Ma} pour une pr\'esentation g\'en\'erale de la th\'eorie des r\'eseaux  dans les groupes de Lie.

\section{Int\'erêt p\'edagogique de la d\'ecomposition effective}

Si on r\'efl\'echit bien, en dehors de ce qui a trait aux op\'erateurs normaux,
il y a deux applications principales d'un cours d'alg\`ebre lin\'eaire
g\'en\'erale au niveau L2/L3.

\smallskip{}

\begin{enumerate}
\item Le calcul des puissances d'une matrice carr\'ee $A$

\smallskip{}

\item Le calcul des exponentielles $e^{tA}$ pour la r\'esolution des syst\`emes
lin\'eaires $Y'(t)=AY(t)+B(t),$ en utilisant la formule

\[
Y(t)=e^{(t-t_{0})A}Y(t_{0})+\int_{t_{0}}^{t}e^{(t-s)A}B(s)ds.\]

\end{enumerate}
La d\'ecomposition de Jordan permet de se ramener
pour ces deux calculs au cas des matrices diagonalisables, puisque
la partie diagonalisable $D$ et la partie nilpotente $N$ de $A$
commutent. On a alors

\[
A^{m}=\sum\limits _{j=0}^{\inf(k,m)}C_{m}^{j}D^{m-j}N^{j},\]

\[
e^{tA}=e^{tD}e^{tN}=e^{tD}\sum_{j=0}^{k}\frac{t^{j}}{j!}N^{j},\]

où $k\ge0$ d\'esigne le plus petit entier tel que $N^{k+1}=0.$

Comme l'unicit\'e de la d\'ecomposition implique que les parties diagonalisables
et nilpotentes d'une matrice \`a coefficients r\'eels sont \`a coefficients
r\'eels, on voit aussi que le fait que les matrices r\'eelles sym\'etriques
sont diagonalisables sur $\mathbb{R}$ d\'ecoule directement du th\'eor\`eme
de d\'ecomposition de Jordan : il suffit de v\'erifier que toute matrice
r\'eelle sym\'etrique et nilpotente est nulle.

Ces consid\'erations, issues de discussions sur le rôle d\'ecomposition
de Jordan dans la th\'eorie des groupes alg\'ebriques pendant la pr\'eparation
de la th\`ese du premier auteur \cite{Cou}, avaient conduit le second
auteur \`a donner \`a la d\'ecomposition de Jordan-Chevalley un rôle central
dans son cours d'alg\`ebre lin\'eaire \`a l'Universit\'e Bordeaux 1 \`a la fin
des ann\'ees 80. Les mêmes raisons ont conduit \`a reprendre depuis 1997
ce point de vue dans le cours destin\'e aux \'el\`eves de premi\`ere ann\'ee
de l'\'ecole d'ing\'enieurs  ESTIA (Ecole Sup\'erieure des Technologies
Industrielles Avanc\'ees), situ\'ee au Pays Basque sur la technopole Izarbel
\`a Bidart, devant des promotions d'\'etudiants d'origines diverses (dont
les effectifs en 1e ann\'ee ont progress\'e de 19 \'etudiants en 1997 \`a
plus de 150 \`a la rentr\'ee 2010) et avec un horaire limit\'e. La m\'ethode
bas\'ee sur l'algorithme de Newton est enseign\'ee depuis 2006, suite
\`a des discussions du troisi\`eme auteur avec des responsables de la
pr\'eparation \`a l'Agr\'egation \`a Bordeaux qui avaient inclus cet algorithme
dans \cite{FrMa}.

La m\'ethode de Newton, sous sa version num\'erique, est connue des \'el\`eves
issus de la fili\`ere MP, mais elle est nouvelle pour la grande majorit\'e
des \'etudiants de 1e ann\'ee de l'ESTIA. Cette m\'ethode permet bien sûr
de faire un lien int\'eressant entre analyse et alg\`ebre, et de
citer dans un cours d'alg\`ebre lin\'eaire H\'eron d'Alexandrie, qui connaissait
l'algorithme

\[
x_{n+1}=\frac{x_{n}^{2}+2}{2x_{n}}\]

qui permet avec $x_{0}=2$ de construire une suite convergeant rapidement
vers $\sqrt{2},$ et n'est autre que l'algorithme obtenu par la m\'ethode
de Newton appliqu\'e \`a la fonction $f:x\to x^{2}-2.$ L'algorithme de
Newton donne aussi une bonne occasion de faire faire aux \'etudiants
l'exercice de programmation simple permettant de faire le calcul effectif
de la d\'ecomposition de Jordan-Chevalley sur d'assez grosses matrices.

La m\'ethode bas\'ee sur le th\'eor\`eme chinois pr\'esente un int\'erêt propre,
car le th\'eor\`eme chinois est un r\'esultat qui passe assez bien aupr\`es
des \'etudiants, et les deux m\'ethodes fonctionnent \'evidement tr\`es bien
en dimension $2$ ou $3.$ Pour le cas $n=2$ le th\'eor\`eme chinois
donne directement $D=\lambda I$ dans le cas d'une matrice $U$ poss\'edant
une valeur propre double, et on a de même $D=\lambda I$ si $U\in{\mathcal{M}}_{n}(k)$
poss\`ede une valeur propre $\lambda$ de multiplicit\'e $n.$ Dans le
cas d'une matrice $3\times3$ poss\'edant une valeur propre simple $\lambda_{1}$
et une valeur propre double $\lambda_{2},$ on peut r\'esoudre imm\'ediatement
le syst\`eme $\left\{ \begin{array}{lr}
p\equiv\lambda_{1}\ \  & \mbox{mod}\ x-\lambda_{1}\\
p\equiv\lambda_{2}\  & \mbox{mod}\ (x-\lambda_{2})^{2}\end{array}\right.$ en posant $p=\lambda_{2}+\alpha(x-\lambda_{2})^{2}.$ La condition
$p(\lambda_{1})=\lambda_{1}$ donne $\alpha=\frac{1}{\lambda_{1}-\lambda_{2}},$
soit $D=\lambda_{2}I+\frac{(U-\lambda_{2})^{2}}{\lambda_{1}-\lambda_{2}}.$
Dans ce cas comme $p_{U}=(x-\lambda_{1})(x-\lambda_{2})^{2},$ la
m\'ethode de Newton donne $\tilde{p}_{U}=(x-\lambda_{1})(x-\lambda_{2}),$
$\overline{p}_{U}=(x-\lambda_{2}).$ Un inverse de $q$ de $p'_{U}=2x-(\lambda_{1}-\lambda_{2})$
modulo $\overline{p}_{U}$ est donc donn\'e par la formule $q=\frac{1}{\lambda_{2}-\lambda_{1}},$
ce qui donne $D=D_{1}=U-\frac{(U-\lambda_{1})(U-\lambda_{2})}{\lambda_{2}-\lambda_{1}}=\lambda_{2}I+\frac{(U-\lambda_{2})^{2}}{\lambda_{1}-\lambda_{2}}.$

\smallskip{}

Une façon classique de construire la d\'ecomposition de Jordan-Chevalley
est d'utiliser la d\'ecomposition spectrale $E=\oplus_{1\le i\le s}Ker(u-\lambda_{i})^{n_{i}}$
associ\'ee \`a un endomorphisme $u$ sur un $k$-espace vectoriel $E$
dont le polynôme caract\'eristique $p_{u}=\prod_{1\le i\le s}(x-\lambda_{i})^{n_{i}}$
est scind\'e sur $k.$ La partie diagonalisable $d$ de $u$ est alors
donn\'ee par la formule $d=\sum_{i=1}\lambda_{j}P_{i},$ où $P_{i}$
est la projection de $E$ sur $Ker(u-\lambda_{i}I_{E})^{n_{i}}$ de
noyau \'egal \`a $\Pi_{j\neq i}Ker(x-\lambda_{i}I_{E})^{n_{j}}.$ Mais
cette formule est exactement celle donn\'ee par le th\'eor\`eme chinois
: pour r\'esoudre le syst\`eme

\[
\left(\mathcal{S}\right):\:\left\{ \begin{array}{c}
h\equiv\lambda_{1}\:\:{\rm mod}\ \left(x-\lambda_{1}\right)^{n_{1}}\\
h\equiv\lambda_{2}\:\:{\rm mod}\ \left(x-\lambda_{2}\right)^{n_{2}}\\
\vdots\\
\vdots\\
h\equiv\lambda_{s}\:\:{\rm mod}\ \left(x-\lambda_{s}\right)^{n_{s}}\end{array}\right.\]
 on construit en utilisant le th\'eor\`eme de Bezout un inverse $E_{i}$
de $(x-\lambda_{i})^{n_{i}}$ modulo $\prod_{j\neq i}(x-\lambda_{j})^{n_{j}},$
et on pose $h=\sum_{i=1}^{s}\lambda_{i}E_{i}(x-\lambda_{i})^{n_{i}}.$
On a donc $d=\sum_{i=1}^{s}\lambda_{i}E_{i}(u)(u-\lambda_{i})^{n_{i}},$
ce qui est la formule pr\'ec\'edente puisque $P_{i}=E_{i}(u)(u-\lambda_{i}I_{E})^{n_{i}}.$

Remarquons enfin que si $p:=(x-\lambda_{1})^{n_{1}}\dots(x-\lambda_{s})^{n_{s}}\in k[x],$
le fait que le syst\`eme ${\mathcal{S}}$ admette une solution $h\in k[x]$
(qui r\'esulte bien sûr de la th\'eorie de Galois, ou du fait qu'une solution
peut être obtenue par l'algorithme de la section 4) est a priori assez
\'evident. On peut remplacer ${\mathcal{S}}$ par ${\mathcal{S}}'$
en remplaçant une solution de ${\mathcal{S}}'$ par le reste de sa
division par $p.$ En posant $n=$max$_{1\le i\le s}n_{j},$ et en
calculant les polynômes $E_{i}$ associ\'es \`a ${\mathcal{S}}'$ par
l'algorithme d'Euclide \'etendu on voit que les coefficients de $h=\sum_{i=1}^{s}\lambda_{i}E_{i}(x-\lambda_{i})^{n}$
sont des quotients de polynômes sym\'etriques de $s$ variables \'evalu\'es
en $(\lambda_{1},\dots,\,\lambda_{s}).$ Ils peuvent donc s'exprimer
comme des fonctions rationnelles des coefficients de $p,$ qui appartiennent
\`a $k.$

En conclusion les auteurs sont convaincus que la d\'ecomposition de
Jordan-Chevalley et son calcul effectif, loin d'être de simples exercices,
doivent recevoir toute l'attention qu'ils m\'eritent et en particulier
jouer un rôle central dans le programme d'alg\`ebre lin\'eaire d'une licence
de Math\'ematiques Fondamentales ou Appliqu\'ees. 

\begin{center}REMERCIEMENTS\end{center}

Les auteurs remercient Pierre de la Harpe, qui a attir\'e leur attention sur les travaux de Venkataramana, et dont les conseils ont permis d'am\'eliorer sensiblement la pr\'esentation de cet article.



\end{document}